\documentclass[]{svmult}

\usepackage{type1cm}
\usepackage{makeidx}
\usepackage{graphicx}
\usepackage{multicol}
\usepackage[bottom]{footmisc}

\usepackage{newtxtext}
\usepackage[varvw]{newtxmath}

\usepackage{amsmath,amsfonts,amssymb,mathtools}

\makeindex

\makeatletter

\DeclareFontFamily{OMX}{MnSymbolE}{}
\DeclareSymbolFont{MnLargeSymbols}{OMX}{MnSymbolE}{m}{n}
\SetSymbolFont{MnLargeSymbols}{bold}{OMX}{MnSymbolE}{b}{n}
\DeclareFontShape{OMX}{MnSymbolE}{m}{n}{
    <-6>  MnSymbolE5
  <6-7>  MnSymbolE6
  <7-8>  MnSymbolE7
  <8-9>  MnSymbolE8
  <9-10> MnSymbolE9
  <10-12> MnSymbolE10
  <12->   MnSymbolE12
}{}
\DeclareFontShape{OMX}{MnSymbolE}{b}{n}{
    <-6>  MnSymbolE-Bold5
  <6-7>  MnSymbolE-Bold6
  <7-8>  MnSymbolE-Bold7
  <8-9>  MnSymbolE-Bold8
  <9-10> MnSymbolE-Bold9
  <10-12> MnSymbolE-Bold10
  <12->   MnSymbolE-Bold12
}{}

\let\llangle\@undefined
\let\rrangle\@undefined
\DeclareMathDelimiter{\llangle}{\mathopen}{MnLargeSymbols}{'164}{MnLargeSymbols}{'164}
\DeclareMathDelimiter{\rrangle}{\mathclose}{MnLargeSymbols}{'171}{MnLargeSymbols}{'171}

\let\llbracket\@undefined
\let\rrbracket\@undefined
\DeclareMathDelimiter{\llbracket}{\mathopen}{MnLargeSymbols}{'102}{MnLargeSymbols}{'102}
\DeclareMathDelimiter{\rrbracket}{\mathclose}{MnLargeSymbols}{'107}{MnLargeSymbols}{'107}

\makeatother

\newcommand{\jump}[1]{\ensuremath{\left\llbracket #1 \right\rrbracket}}
\newcommand{\mean}[1]{\ensuremath{\left\llangle #1 \right\rrangle}}
\newcommand{\scp}[2]{\ensuremath{\big\langle#1,#2\big\rangle}}

\begin{document}

\title*{Discontinuous Galerkin method for incompressible two-phase flows}
\author{Janick Gerstenberger, Samuel Burbulla and Dietmar Kr\"oner}
\institute{
  Janick Gerstenberger \at
  AAM, Albert-Ludwigs-Universit\"at Freiburg,
  Hermann-Herder-Str. 10, 79104 Freiburg, Germany,\\
  \email{janick.gerstenberger@mathematik.uni-freiburg.de}
  \and
  Samuel Burbulla \at
  IANS, Universit\"at Stuttgart, Pfaffenwaldring 57, 70569 Stuttgart, Germany,\\
  \email{samuel.burbulla@mathematik.uni-stuttgart.de}
  \and
  Dietmar Kr\"oner \at
  AAM, Albert-Ludwigs-Universit\"at Freiburg,
  Hermann-Herder-Str. 10, 79104 Freiburg, Germany,\\
  \email{dietmar@mathematik.uni-freiburg.de}}

\maketitle

\abstract{
  In this contribution we present a local discontinuous Galerkin (LDG) pressure-correction scheme for the incompressible Navier-Stokes equations.
  The scheme does not need penalty parameters and satisfies the discrete continuity equation exactly.
  The scheme is especially suitable for two-phase flow when used with a piecewise-linear interface construction (PLIC) volume-of-fluid (VoF) method and cut-cell quadratures.
  \keywords{Two-phase flow, Incompr. Navier-Stokes Eqs., Discontinuous Galerkin\\[5pt]\textbf{MSC}(2010)\textbf{:} 76T10, 76D45, 76D05, 65M60}
}

\section{Introduction}

Sharp interface models for incompressible two-phase flows have gained in popularity in recent years.
These models combine incompressible flows in the bulk domains with jump conditions along the interfaces that separate the fluids,
which model fluid interactions and surface effects, like surface tension.

Discontinuous Galerkin methods are a popular choice for solving incompressible flow problems due to their local mass conservation property and potentially high order of convergence.

Designing schemes for solving such two-phase flows presents several challenges.
The choice of method for the phase transport has implications for the interface representation and conservation properties.
Volume-of-fluid~\cite{rider98} methods are conservative but their interface representations is in general non-continuous,
while level-set methods have continuous interface representations, but are not conservative by themselves.
There are several further methods that combine features from or generalize volume-of-fluid and level-set methods, but these are comparatively more complex.
Incompressible flow problems have a saddle-point structure, which can make them computationally difficult/expensive to solve.
Splitting methods like the various projection methods where introduced to decouple the problems into simpler ones.
Projections methods replace the saddle-point problem with a advection-diffusion equation for the velocity and Poisson problem for the pressure,
which are computationally simpler to solve.
But for most discontinuous Galerkin discretizations the discrete continuity equation is not satisfied without some postprocessing techniques like the \(H(div)\) reconstruction presented in \cite{piatkowski2018}.
Stability in regards to high coefficient ratios and strong surface effects are further issues that need to be addressed.

We present a discontinuous Galerkin pressure-correction method for incompressible two-phase flow that is robust in regards to coefficient ratios.
The scheme is simple to implement and has shown good results for benchmarks problems and some numerical experiments with realistic data.

This contribution is structured as follows: We first briefly present the model for incompressible two-phase flow without phase transitions.
Next we give some notation, present the modified LDG method without penalization introduced in \cite{neilan2016}
and present the new LDG pressure-correction scheme for incompressible two-phase flows based on it.
Then we present numerical experiments for benchmark problems from the literature~\cite{hysing09}.
Lastly we give some concluding remarks and an outlook on further work.

\section{Model}
The model we consider is a simplification of the sharp interface model presented in \cite{pruess14} without phase transitions.
Let \(\Omega\subset\mathbb{R}^n\), \(n=2, 3\) be a bounded domain, which is divided into two disjunct phases: at time \(t\), phase \(i = 1, 2\) occupies subdomain \(\Omega_i(t) \subset \Omega\).
The boundary between the phase is the (phase) interface \(\Gamma(t) \coloneq \partial\Omega_1\cap\partial\Omega_2\).

Let \(\vec{u}\) denote the velocity field, \(p\) the pressure field, \(\nu\) the outer normal on \(\Omega_1\), \(\kappa\) the mean curvature of \(\Gamma(t)\) and \(\jump{q} \coloneq q_2 -q_1\) the jump of the variable \(q\) across \(\Gamma(t)\).
In addition the constants \(\varrho_1, \varrho_2 > 0\) denote the densities, \(\mu_1, \mu_2 > 0\) the viscosities of the phases, \(\vec{g}\) the gravitational acceleration and \(\sigma\) the coefficient of surface tension.
In the following we drop the indices on density and viscosity and keep in mind that these coefficients depend on the phase.

Then the model is given by the incompressible Navier-Stokes equations in each phase,
\begin{equation}
  \left.
\begin{aligned}
  \partial_{t} (\varrho \vec{u}) + \nabla\cdot\!{(\varrho \vec{u} \otimes \vec{u} + \tens{T})} &= \varrho\vec{g}\\
  \nabla\cdot\!{\vec{u}} &= 0\\
\end{aligned}
  \quad\right\}\quad \text{in } \Omega_1, \Omega_2\,,
\end{equation}
with the stress tensor \(\tens{T} = p\tens{I} - 2\umu\,\mathcal{D}{\vec{u}} \), \(\mathcal{D}{\vec{u}} = S(\nabla{\vec{u}}) \coloneq (\nabla{\vec{u}} + \nabla{\vec{u}}^\top) / 2\),
augmented with the following jump conditions
  \begin{equation}
  \left.
\begin{aligned}
  \jump{\tens{T}}\vec{\nu} &= -\sigma \kappa \vec{\nu}\\
  \jump{\vec{u}} &= 0\\
  \end{aligned}
  \quad\right\}\quad \text{on } \Gamma.
\end{equation}
Additionally we impose either no-slip or free-slip conditions on the boundary of \(\Omega\).

\section{Discretization}

To discretize the model above we employ a primal LDG method~\cite{cockburn05}.
It can be formulated in terms of a discrete gradient operator that is composed of the elementwise gradient and a lifting of the jumps into the piecewise discrete space.
By constructing the liftings one order higher than the used discrete space the method is rendered stable without penalty parameters \cite{neilan2016}.
From these building blocks we construct an incremental pressure-projection scheme that satisfies the discrete continuity equation.

Since we require strict mass-conservation we choose to use a PLIC-VoF method~\cite{rider98,diot2016}.
These methods have the disadvantage that the reconstructed interface is in general not continuous and its curvature needs further approximations \cite{cummins05}, but is mass conservative.

We chose this LDG method because with the discrete gradients and no penalty terms  no (explicit) evaluations of integrals over element boundaries are needed.
In our experience the non-continuous discrete interface reconstructions on the inter-element boundaries lead to issues around the interface.

\subsection{Notation and Liftings}
  To derive the discretization we first introduce some notation.
  Let \(\mathcal{T}_h \) be a triangulation of \(\Omega\) into elements \(E\).
  By \(\Sigma_h^I\) we denote the set of all interior intersections \(e\) of elements \(E^-, E^+ \in \mathcal{T}_h\) with \(e = E^- \cap E^+ \neq \emptyset\),
  by \(\Sigma^D_h\) the set of all intersections with Dirichlet boundary values,
  by \(\Sigma^N_h\) the set of intersections with Neumann boundary values
  and by \(\Sigma_h = \Sigma_h^I \cup \Sigma_h^D \cup \Sigma_h^N\) the set of all intersections.

  For \(e \in \Sigma_h^I\), with \(e = E^- \cap E^+\),  we introduce the jump and the average
  \begin{subequations}
    \begin{align}
      \jump{u} &= u_{ |_{E^-}} - u_{ |_{E^+}}\,,\\
      \mean{u} &= \tfrac{1}{2}\left(u_{ |_{E^-}} + u_{ |_{E^+}}\right)
    \end{align}
  \end{subequations}

  We define the discrete spaces of piecewise polynomials of degree \(\leq k\)
  \begin{equation}
    \mathcal{V}_k^d = \left\{\vec{v}\in[L^{2}(\Omega)]^{d}\mid \vec{v}_{ |_E}\in[\mathcal{P}_k(E)]^{d},\,E\in \mathcal{T}_h\right\}
  \end{equation}
  and \(\mathcal{V}_k = \mathcal{V}_k^1\) in the scalar case.
  The piecewise \(L^{2}\) inner product over \(\mathcal{T}_h\) is given by
  \begin{equation}
    \scp{\vec{v}}{\vec{w}} = \sum_{E\in\mathcal{T}_h}\!\! \scp{\vec{v}}{\vec{w}}_E = \sum_{E\in\mathcal{T}_h} \int_E\! \vec{v}\cdot\vec{w}
  \end{equation}

  The gradient lifting operators \(R: \mathcal{V}_k \mapsto \mathcal{V}_{k+1}^{d}\) and \(R_a: \mathcal{V}_k \mapsto \mathcal{V}_{k+1}^{d}\)
   are defined by
  \begin{subequations}
    \begin{equation}
      \scp{R(\jump{v})}{\vec{w}}
      = \sum_{e \in \Sigma_h^I} \int_{e}\!\jump{v}\mean{\vec{w}}\cdot\vec{n}_e\,,
    \end{equation}
    \begin{equation}
      \scp{R_a(\jump{v})}{\vec{w}}
      = \scp{R(\jump{v})}{\vec{w}}
      + \sum_{e \in \Sigma_h^D} \int_{e}\!(v-a) \vec{w}\cdot\vec{n}_e\,,
    \end{equation}
  \end{subequations}
  for all \(\vec{w}\in\mathcal{V}_{k+1}^d\).
  Similarly, the divergence lifting operators \(M: \mathcal{V}_{k+1}^d \mapsto \mathcal{V}_{k}\), \(M_{\vec{b}}: \mathcal{V}_{k+1}^d \mapsto \mathcal{V}_{k}\)
  are defined by
  \begin{subequations}
    \begin{equation}
      \scp{M(\jump{\vec{v}})}{w}
      = \sum_{e \in \Sigma_h^I} \int_{e}\!\jump{\vec{v}}\cdot\vec{n}_e\,\mean{w}\,,
    \end{equation}
    \begin{equation}
      \scp{M_{\vec{b}}(\jump{\vec{v}})}{w}
      = \scp{M(\jump{\vec{v}})}{w}
      + \sum_{e \in \Sigma_h^D} \int_{e}\!(\vec{v}-\vec{b})\cdot\vec{n}_e\,w\,,
    \end{equation}
  \end{subequations}
  for all \(w\in\mathcal{V}_k\). Here \(a, \vec{b}\) are the respective Dirichlet boundary conditions.

  With the liftings we can now define the lifted DG gradient and divergence
  \begin{subequations}
    \begin{align}
      \overline{\nabla}v &= \nabla_{\!h} v - R(\jump{v}), &\overline{\nabla}_{\! g}v &= \nabla_{\!h} v - R_g(\jump{v}), & v &\in\mathcal{V}_k,\\
      \overline{\nabla}\!\cdot\!\vec{v} &= \nabla_{\! h}\!\cdot\!\vec{v} - M(\jump{\vec{v}}), & \overline{\nabla}_{\!g}\!\cdot\!\vec{v} &= \nabla_{\! h}\!\cdot\!\vec{v} - M_g(\jump{\vec{v}}), & \vec{v} &\in\mathcal{V}_k^d.
    \end{align}
  \end{subequations}
  The lifted derivatives with homogeneous Dirichlet boundary (\(a =0, \vec{b} = 0\)) satisfy the following discrete integration-by-parts identities, as shown in \cite{feng2016}:
  \begin{subequations}
    \begin{align}
      \scp{\overline{\nabla}_{\!0}\!\cdot\!\vec{v}}{w} &= -\scp{\vec{v}}{\overline{\nabla}w}\label{eq:part-int-a}\\
      \scp{\overline{\nabla}\!\cdot\!\vec{v}}{w} &= -\scp{\vec{v}}{\overline{\nabla}_{\!0}w}\label{eq:part-int-b}
    \end{align}
  \end{subequations}
  for  all \(\vec{v} \in \mathcal{V}_{k+1}^d\), \(w\in \mathcal{V}_k\).
  This means the lifted derivatives are adjoint to each other.
  These identies are useful for defining of projection methods with respect to the lifted derivatives.

\subsection{Unpenalized LDG Scheme}

  As a simple example we now consider the discretization of the Poisson equation \(-\Delta{u} = f\) with homogeneous Dirichlet boundary data.
  The modified LDG~method~\cite{neilan2016} reads: Find \(u\in\mathcal{V}_k\) such that
  \begin{equation}
    \scp{\overline{\nabla}_{\!0}u}{\overline{\nabla}_{\!0}v} = \scp{f}{v}\qquad \forall v\in\mathcal{V}_k.
  \end{equation}
  This scheme with order \(k+1\) liftings is stable without adding penalty terms and can also be written in a ``strong''-form by using the integration-by-parts identity~\eqref{eq:part-int-b}
  \begin{equation}
    -\scp{\underline{\Delta}_{\mkern1.5mu 0} u}{v} \coloneq -\scp{\overline{\nabla}\!\cdot\!\overline{\nabla}_{\!0}u}{v} = \scp{f}{v}\qquad \forall v\in\mathcal{V}_k.
  \end{equation}

\subsection{Two-Phase LDG Scheme}

We now present our primal LDG pressure-correction scheme for two-phase flow.

The scheme first split into an explicit (linearized) advection step for the momentum/velocity and the phase transport and an implicit Stokes step.
Because the momentum and the phase interface are transported at the same rate, we can formulate the explicit step in the velocity,
which also means this part of the scheme does not depend on the phase interface.
The Stokes step is then further split into an implicit momentum step, a pressure Poisson step and an update step.

To sharply resolve the phases the terms containing a phase dependent coefficient are integrated using cut-cell quadratures.
The cut-cell quadratures are constructed by cutting the elements containing an interface reconstruction with its interface reconstruction, subtriangulating the part of each phase and using standard simplex quadratures in the the subtriangulations.
Surface tension effect are included by integration of the jump condition over segments of the interface reconstruction.
Using the lifted derivatives we eliminate (explicit) evaluations of integrals containing phase dependent coefficients along element boundaries, which in our experience can cause instabilities.

To simplify the presentation we restrict ourself to first order time stepping (IMEX Euler) and assume the phase transport/interface reconstruction is given:
Let \(\Gamma_h^{ n+1}\) denote the set of all interface reconstruction at time \(t^{n+1}\) and \(\kappa^{n+1}\) the approximated interface mean curvature at time \(t^{n+1}\).

The complete scheme for the velocity and pressure at time \(t^{n+1}\) then reads as:
Find \(\vec{u}^{n+1}\in\mathcal{V}_{k+1}^d\), \(p^{n+1}\in\mathcal{V}_{k}\) such that
\begin{subequations}
  \begin{equation}
    \frac{1}{\Delta t}\scp{\hat{\vec{u}}^*-\vec{u}^n}{\vec{v}}
    - \scp{\vec{u}^n\otimes\vec{w}^n}{\nabla{\vec{v}}}
    + \!\sum_{e\in\Sigma_h}\!\int_{e}\!\!F_e(\vec{u}^n;\vec{w}^n)\cdot\vec{v}
    = \scp{\vec{g}}{\vec{v}}\,,
  \end{equation}
  \begin{align}
    \frac{1}{\Delta t}\scp{\varrho(\vec{u}^{*}-\hat{\vec{u}}^*)}{\vec{v}}
    - \scp{2\umu\, S(\overline{\nabla}_{\!0}\vec{u}^{*})}{\overline{\nabla}_{\!0}\vec{v}}
    &= -\scp{\overline{\nabla}p^n}{\vec{v}}
    + \!\!\!\sum_{\gamma\in\Gamma_h^{n+1}}\!\int_{\gamma}\!\!\sigma \kappa^{n+1} \vec{v}\!\cdot\!\vec{n}_\gamma\,,\label{eq:implicit-a}\\
    \scp{\frac{1}{\varrho}\overline{\nabla}p^*}{\overline{\nabla}q} &= \frac{1}{\Delta t}\scp{\overline{\nabla}_{\!0}\!\cdot\!\vec{u}^{*}}{q}\,,\label{eq:implicit-b}\\
    \vec{u}^{n+1} = \vec{u}^{*} + \frac{\Delta t}{\varrho}\overline{\nabla}p^*\,&, \qquad p^{n+1} = p^{n} + p^*\,,\label{eq:update}
  \end{align}
\end{subequations}
for all \(\vec{v}\in \mathcal{V}_{k+1}^d\), \(q\in\mathcal{V}_k\).
Here \(F_e\) is a suitable numerical flux (e.g.\ local Lax-Friedrichs flux), \(\vec{w}^n = \vec{u}^n\) is the transport velocity and all terms in eqs.~\eqref{eq:implicit-a}, \eqref{eq:implicit-b} containing \(\varrho\), \(\umu\) are integrated with cut-cell quadratures with respect to \(\Gamma_h^{n+1}\).

We note that here the lifted gradient/divergence maps into the dg space with increased/decreased polynomial order of its argument function space (e.g. \(\overline{\nabla}_{\!0}\vec{u}\in\mathcal{V}_{k+2}^{d\times d}\), \(\overline{\nabla}_{\!0}\!\cdot\!\vec{u}\in\mathcal{V}_{k}, \)).

The resulting discrete velocity field is exactly divergence-free with respect to the discrete lifted divergence operator for arbitrary density ratios (up to the accurracy of the solver used).
Applying the discrete DG divergence to the velocity field update and using eq.~\eqref{eq:part-int-a}  recovers eq.~\eqref{eq:implicit-b}.

\section{Numerical experiments}

The scheme has been implemented in DUNE-FEM~\cite{dedner10}, which is part of the DUNE (Distributed and Unified Numerics Environment) framework~\cite{dune24:16}.
The discretization used here are orthonormal \(\mathcal{P}_2\) (resp.\ \(\mathcal{P}_1\)) dg elements for the velocity (resp.\ pressure) on a cartesian grid for space and a BDF2 scheme for time stepping.
The implicit substeps are solved with a GMRES preconditioned with either an incomplete \(LU\)- or \(LDL^\mathsf{T}\)-decomposition.
The phase transport is solved using a VoF method with a geometric flux, also implemented in DUNE.

We consider the benchmark problems presented in \cite{hysing09}.
The general setup is the following:
The computational domain is given by \( \Omega = [0,1]\times[0,2] \subset \mathbb{R}^2 \), \(\Omega_2\) is a circular subdomain with diameter \(d = 0.5\) centered at \((0.5, 0.5)^\top\) and \(\Omega_1 \coloneq \Omega \setminus \Omega_2\).
The no-slip boundary condition is prescribed at the top and bottom boundaries, the free-slip condition is prescribed on the vertical boundaries and the system is initially at reset.

\begin{figure}[tb]
    \centering
      \includegraphics[trim=0 10cm 0 4cm, clip, width=0.4\linewidth]{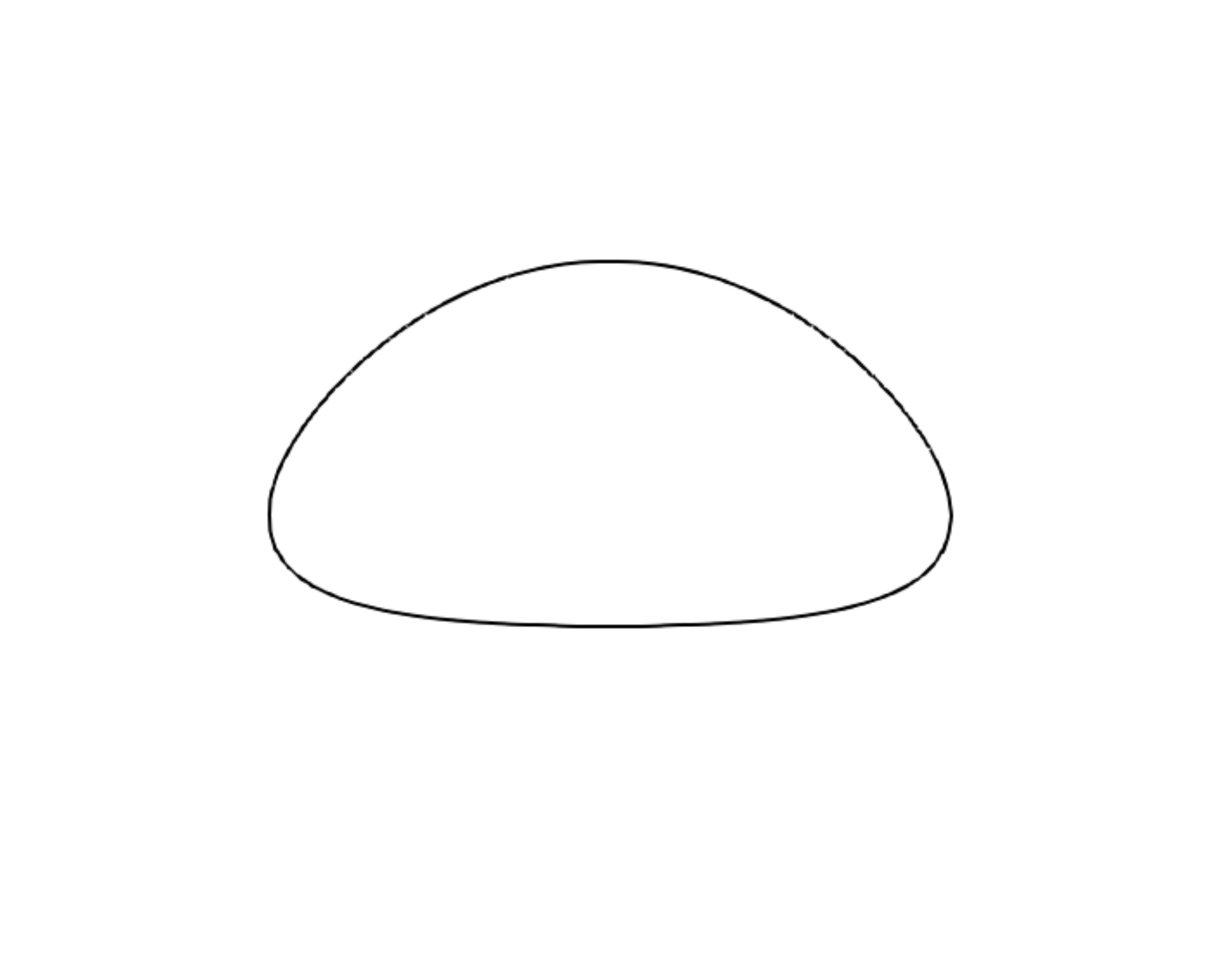}
      \includegraphics[trim=0 10cm 0 4cm, clip, width=0.4\linewidth]{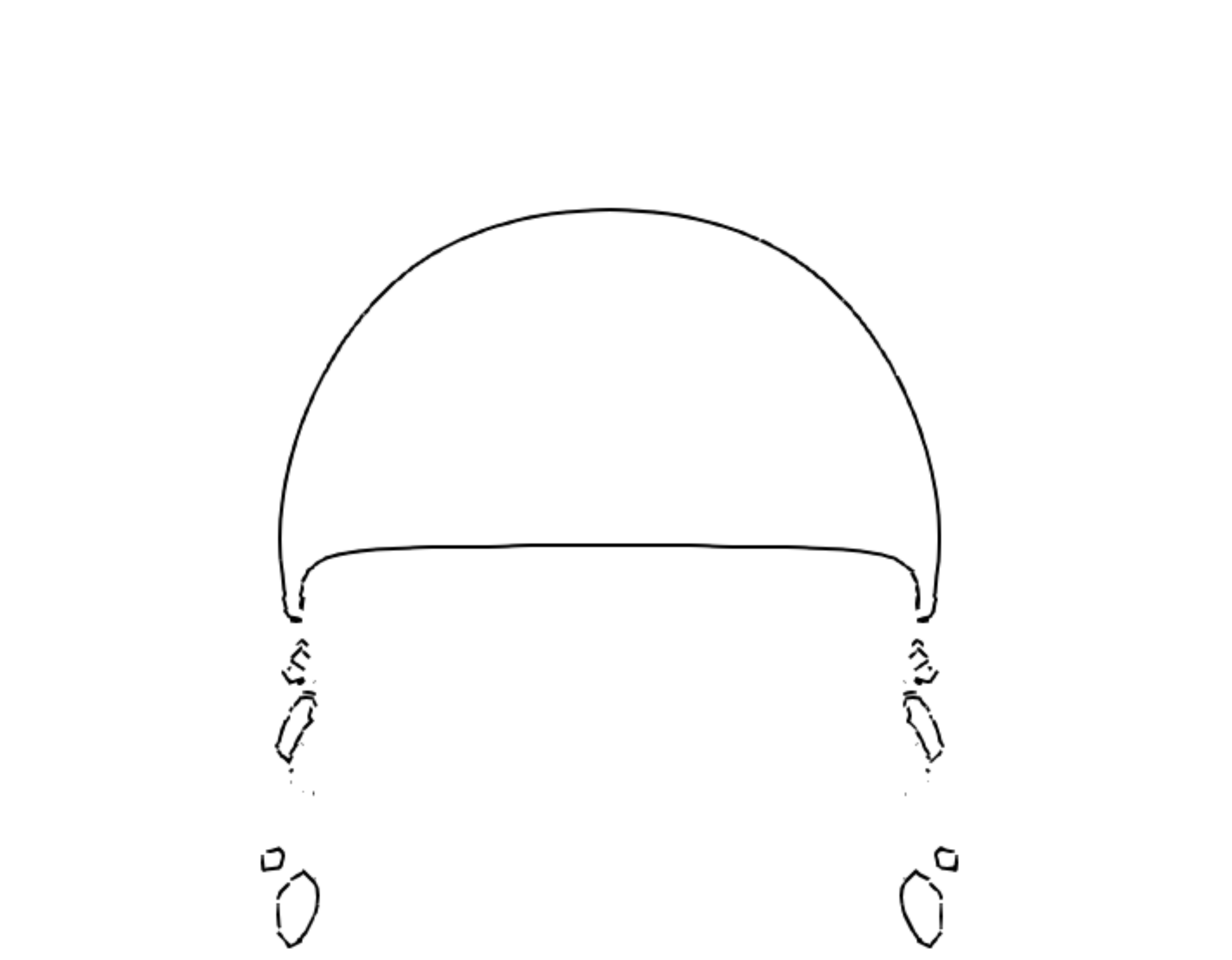}
    \caption{Bubble shape of test case 1 (left) and test case 2 (right) at the final time \(t=3\) on a cartesian mesh with 80 by 160 cells}%
    \label{fig:1}
\end{figure}

The fluid coefficients for the test problems are given in table~\ref{tab:1}.
Test case 1 results in a ellipsoidal bubble, while in test case 2 the bubble gets deformed significantly and eventually break ups occur.
Our results are in good agreement for the shape (see Fig.~\ref{fig:1}), the rise velocity and the center of mass of the bubble with the results presented in \cite{hysing09}.
The `circularity' can not be compared directly since our interface representation makes measuring the perimeter of the bubble difficult, nonetheless it is in general agreement.

\begin{table}[bt]
  \caption[Physical parameters and dimensionless numbers defining the test cases.]{%
    Physical parameters and dimensionless numbers defining the test cases.\\
    \(Re = \varrho_1 U_g d / \umu_1\), \(Eo = \varrho_1 U_g^2 d / \sigma \), \(U_g = \sqrt{gd}\).}%
    \label{tab:1}
  \centering
  \begin{tabular*}{\textwidth}{@{\extracolsep{\fill} } l*{10}{r}}
    \svhline\noalign{\smallskip}
    case & \(\varrho_1\) &\(\varrho_2\) & \(\umu_1\) &\(\umu_2\) &\(g\) &\(\sigma\) & \(Re\) &\(Eo\) & \(\varrho_1 / \varrho_2\) & \(\umu_1 / \umu_2\) \\
    \noalign{\smallskip}\hline\noalign{\smallskip}
    1 & \(1000\) & \(100\) & \(10\) & \(1\) & \(0.98\) & \(24.5\hphantom{6}\) & \(35\) & \(10\) & \(10\) & \(10\)\\
    2 & \(1000\) & \(1\) & \(10\) & \(0.1\) & \(0.98\) & \(1.96\) & \(35\) & \(125\) & \(1000\) & \(100\)\\
    \noalign{\smallskip}\svhline
  \end{tabular*}
\end{table}

\begin{figure}[htb]
  \centering
    \includegraphics[width=0.75\linewidth]{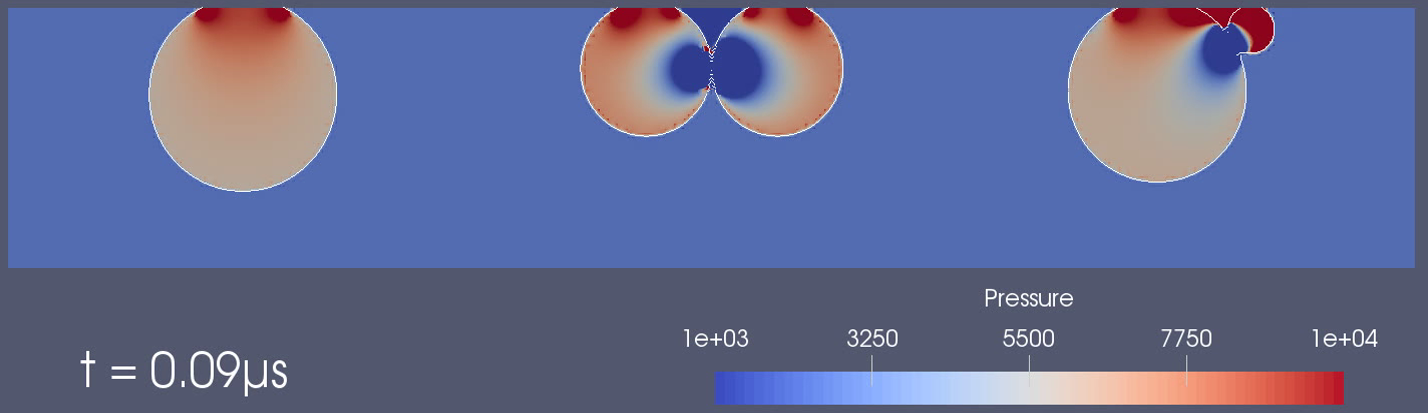}
    \includegraphics[width=0.75\linewidth]{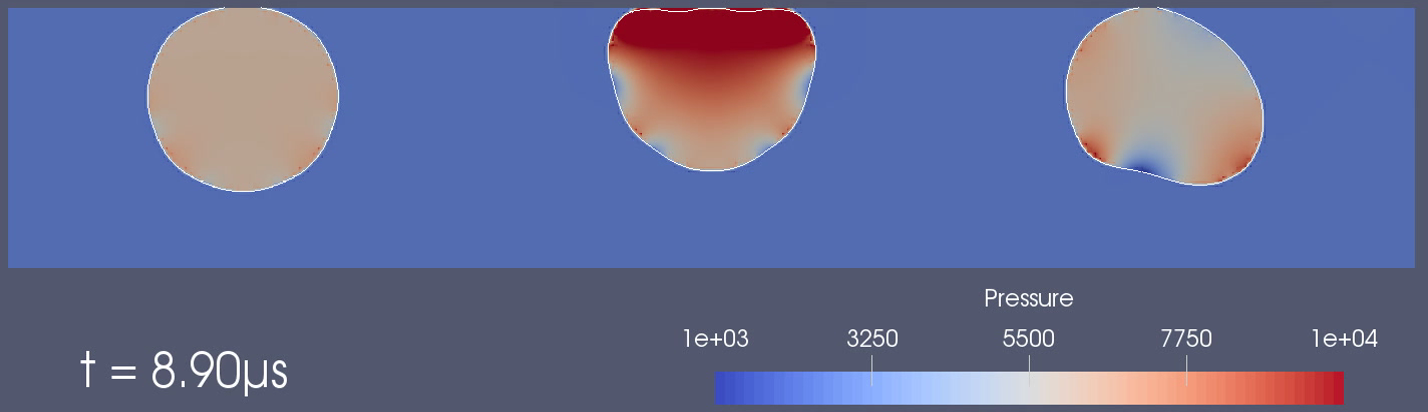}
    \includegraphics[width=0.75\linewidth]{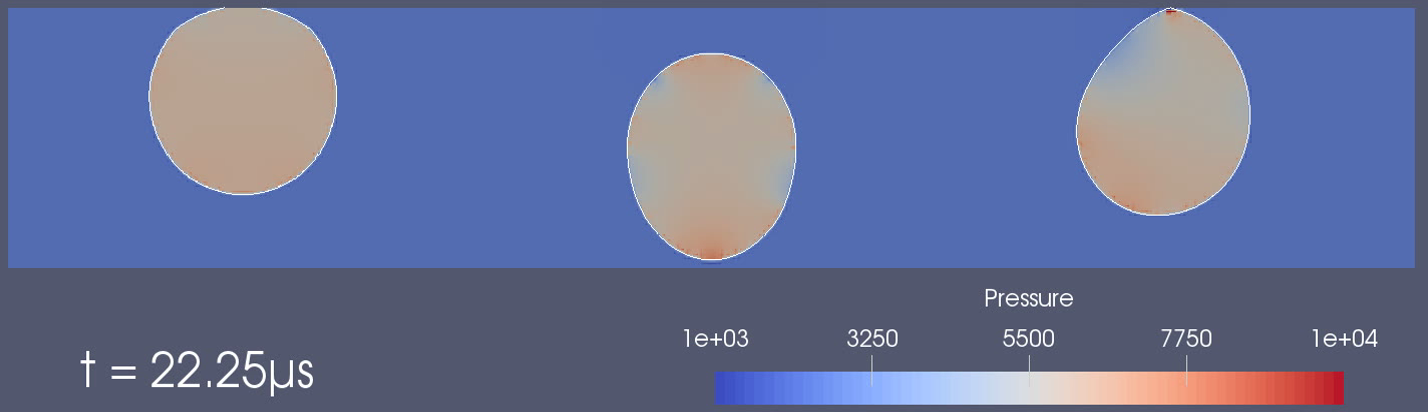}
  \caption{Pressure profile of water droplets in a steam atmosphere on super-hydrophobic surface (contact angle \(\alpha = 165^\circ\)) in free fall on a cartesian mesh with 300 by 100 cells, \(h = 10^{-4} \mathrm{m}\).}%
  \label{fig:2}
\end{figure}

Figure~\ref{fig:2} shows an experiment of the effect of droplet merging on a super-hydrophobic surface in a water-steam system in free-fall.
The initial configuaration consist of two equal-sized droplets, two unequal-sized droplets and a single droplet as reference.
When droplets of similar size merge the force of surface tension is strong enough that the droplet ``jumps''.

\section{Conclusions and Outlook}

In this contribution we presented a primal local discontinuous Galerkin pressure correction scheme suitable for two-phase flows with high density and viscosity ratios.
The interface is sharply resolved by using cut-cell quadratures.
Numerical experiments show that the scheme agrees with other results presented in the literature.

Further work includes extending the model and the scheme to include the effects of phase transition and investigating if other, more compact, DG methods can recast in the same form as the presented LDG method.

\begin{acknowledgement}
  The first author would like to thank Tobias Malkmus for providing the initial code base for the fluid solver and explaining its subtleties, Martin Nolte for his ideas and help in modifying said code.
\end{acknowledgement}

\bibliographystyle{spmpsci}
\bibliography{references}

\end{document}